
\font \fivesans               = cmss10 at 5pt

\font \ninebf                 = cmbx9
\font \ninei                  = cmmi9
\font \nineit                 = cmti9
\font \ninerm                 = cmr9
\font \ninesans               = cmss10 at 9pt
\font \ninesl                 = cmsl9
\font \ninesy                 = cmsy9
\font \ninett                 = cmtt9
\font \sevensans              = cmss10 at 7pt
\font \sixbf                  = cmbx6
\font \sixi                   = cmmi6
\font \sixrm                  = cmr6
\font \sixsans                = cmss10 at 6pt
\font \sixsy                  = cmsy6
\font \smallescriptfont       = cmr5 at 7pt
\font \smallescriptscriptfont = cmr5
\font \smalletextfont         = cmr5 at 10pt

\font \tafonts                = cmbx7  scaled\magstep2
\font \tafontss               = cmbx5  scaled\magstep2
\font \tafontt                = cmbx10 scaled\magstep2
\font \tams                   = cmmib10
\font \tamss                  = cmmib10 scaled 700
\font \tamt                   = cmmib10 scaled\magstep2
\font \tass                   = cmsy7  scaled\magstep2
\font \tasss                  = cmsy5  scaled\magstep2
\font \tast                   = cmsy10 scaled\magstep2

\font \tbfonts                = cmbx7  scaled\magstep1
\font \tbfontss               = cmbx5  scaled\magstep1
\font \tbfontt                = cmbx10 scaled\magstep1
\font \tbms                   = cmmib10 scaled 833
\font \tbmss                  = cmmib10 scaled 600
\font \tbmt                   = cmmib10 scaled\magstep1
\font \tbss                   = cmsy7  scaled\magstep1
\font \tbsss                  = cmsy5  scaled\magstep1
\font \tbst                   = cmsy10 scaled\magstep1

\font \tensans                = cmss10

\vsize=19.5cm
\hsize=12.1cm
\hfuzz=2pt
\tolerance=500
\abovedisplayskip=3 mm plus6pt minus 4pt
\belowdisplayskip=3 mm plus6pt minus 4pt
\abovedisplayshortskip=0mm plus6pt minus 2pt
\belowdisplayshortskip=2 mm plus4pt minus 4pt
\predisplaypenalty=0
\clubpenalty=10000
\widowpenalty=10000
\frenchspacing
\newdimen\oldparindent\oldparindent=1.5em
\parindent=1.5em
\skewchar\ninei='177 \skewchar\sixi='177
\skewchar\ninesy='60 \skewchar\sixsy='60
\hyphenchar\ninett=-1
\def\newline{\hfil\break}%
\catcode`@=11
\def\folio{\ifnum\pageno<\z@
\uppercase\expandafter{\romannumeral-\pageno}%
\else\number\pageno \fi}
\catcode`@=12 
  \mathchardef\Gamma="0100
  \mathchardef\Delta="0101
  \mathchardef\Theta="0102
  \mathchardef\Lambda="0103
  \mathchardef\Xi="0104
  \mathchardef\Pi="0105
  \mathchardef\Sigma="0106
  \mathchardef\Upsilon="0107
  \mathchardef\Phi="0108
  \mathchardef\Psi="0109
  \mathchardef\Omega="010A

\def\sq{\hbox{\rlap{$\sqcap$}$\sqcup$}}

\def\utw{\smash{\rlap{\lower5pt\hbox{$\sim$}}}}
\def\udtw{\smash{\rlap{\lower6pt\hbox{$\approx$}}}}

\def\diameter{{\ifmmode\mathchoice
{\ooalign{\hfil\hbox{$\displaystyle/$}\hfil\crcr
{\hbox{$\displaystyle\mathchar"20D$}}}}
{\ooalign{\hfil\hbox{$\textstyle/$}\hfil\crcr
{\hbox{$\textstyle\mathchar"20D$}}}}
{\ooalign{\hfil\hbox{$\scriptstyle/$}\hfil\crcr
{\hbox{$\scriptstyle\mathchar"20D$}}}}
{\ooalign{\hfil\hbox{$\scriptscriptstyle/$}\hfil\crcr
{\hbox{$\scriptscriptstyle\mathchar"20D$}}}}
\else{\ooalign{\hfil/\hfil\crcr\mathhexbox20D}}%
\fi}}


\def\bbbc{{\mathchoice {\setbox0=\hbox{$\displaystyle\rm C$}\hbox{\hbox
to0pt{\kern0.4\wd0\vrule height0.9\ht0\hss}\box0}}
{\setbox0=\hbox{$\textstyle\rm C$}\hbox{\hbox
to0pt{\kern0.4\wd0\vrule height0.9\ht0\hss}\box0}}
{\setbox0=\hbox{$\scriptstyle\rm C$}\hbox{\hbox
to0pt{\kern0.4\wd0\vrule height0.9\ht0\hss}\box0}}
{\setbox0=\hbox{$\scriptscriptstyle\rm C$}\hbox{\hbox
to0pt{\kern0.4\wd0\vrule height0.9\ht0\hss}\box0}}}}
\def\bbbe{{\mathchoice {\setbox0=\hbox{\smalletextfont e}\hbox{\raise
0.1\ht0\hbox to0pt{\kern0.4\wd0\vrule width0.3pt
height0.7\ht0\hss}\box0}}
{\setbox0=\hbox{\smalletextfont e}\hbox{\raise
0.1\ht0\hbox to0pt{\kern0.4\wd0\vrule width0.3pt
height0.7\ht0\hss}\box0}}
{\setbox0=\hbox{\smallescriptfont e}\hbox{\raise
0.1\ht0\hbox to0pt{\kern0.5\wd0\vrule width0.2pt
height0.7\ht0\hss}\box0}}
{\setbox0=\hbox{\smallescriptscriptfont e}\hbox{\raise
0.1\ht0\hbox to0pt{\kern0.4\wd0\vrule width0.2pt
height0.7\ht0\hss}\box0}}}}
\def\bbbq{{\mathchoice {\setbox0=\hbox{$\displaystyle\rm
Q$}\hbox{\raise
0.15\ht0\hbox to0pt{\kern0.4\wd0\vrule height0.8\ht0\hss}\box0}}
{\setbox0=\hbox{$\textstyle\rm Q$}\hbox{\raise
0.15\ht0\hbox to0pt{\kern0.4\wd0\vrule height0.8\ht0\hss}\box0}}
{\setbox0=\hbox{$\scriptstyle\rm Q$}\hbox{\raise
0.15\ht0\hbox to0pt{\kern0.4\wd0\vrule height0.7\ht0\hss}\box0}}
{\setbox0=\hbox{$\scriptscriptstyle\rm Q$}\hbox{\raise
0.15\ht0\hbox to0pt{\kern0.4\wd0\vrule height0.7\ht0\hss}\box0}}}}
\def\bbbt{{\mathchoice {\setbox0=\hbox{$\displaystyle\rm
T$}\hbox{\hbox to0pt{\kern0.3\wd0\vrule height0.9\ht0\hss}\box0}}
{\setbox0=\hbox{$\textstyle\rm T$}\hbox{\hbox
to0pt{\kern0.3\wd0\vrule height0.9\ht0\hss}\box0}}
{\setbox0=\hbox{$\scriptstyle\rm T$}\hbox{\hbox
to0pt{\kern0.3\wd0\vrule height0.9\ht0\hss}\box0}}
{\setbox0=\hbox{$\scriptscriptstyle\rm T$}\hbox{\hbox
to0pt{\kern0.3\wd0\vrule height0.9\ht0\hss}\box0}}}}
\def\bbbs{{\mathchoice
{\setbox0=\hbox{$\displaystyle     \rm S$}\hbox{\raise0.5\ht0\hbox
to0pt{\kern0.35\wd0\vrule height0.45\ht0\hss}\hbox
to0pt{\kern0.55\wd0\vrule height0.5\ht0\hss}\box0}}
{\setbox0=\hbox{$\textstyle        \rm S$}\hbox{\raise0.5\ht0\hbox
to0pt{\kern0.35\wd0\vrule height0.45\ht0\hss}\hbox
to0pt{\kern0.55\wd0\vrule height0.5\ht0\hss}\box0}}
{\setbox0=\hbox{$\scriptstyle      \rm S$}\hbox{\raise0.5\ht0\hbox
to0pt{\kern0.35\wd0\vrule height0.45\ht0\hss}\raise0.05\ht0\hbox
to0pt{\kern0.5\wd0\vrule height0.45\ht0\hss}\box0}}
{\setbox0=\hbox{$\scriptscriptstyle\rm S$}\hbox{\raise0.5\ht0\hbox
to0pt{\kern0.4\wd0\vrule height0.45\ht0\hss}\raise0.05\ht0\hbox
to0pt{\kern0.55\wd0\vrule height0.45\ht0\hss}\box0}}}}
\def\bbbz{{\mathchoice {\hbox{$\sans\textstyle Z\kern-0.4em Z$}}
{\hbox{$\sans\textstyle Z\kern-0.4em Z$}}
{\hbox{$\sans\scriptstyle Z\kern-0.3em Z$}}
{\hbox{$\sans\scriptscriptstyle Z\kern-0.2em Z$}}}}
\def\qed{\ifmmode\sq\else{\unskip\nobreak\hfil
\penalty50\hskip1em\null\nobreak\hfil\sq
\parfillskip=0pt\finalhyphendemerits=0\endgraf}\fi}
\newfam\sansfam
\textfont\sansfam=\tensans\scriptfont\sansfam=\sevensans
\scriptscriptfont\sansfam=\fivesans
\def\sans{\fam\sansfam\tensans}
\def\stackfigbox{\if
Y\FIG\global\setbox\figbox=\vbox{\unvbox\figbox\box1}%
\else\global\setbox\figbox=\vbox{\box1}\global\let\FIG=Y\fi}
\def\placefigure{\dimen0=\ht1\advance\dimen0by\dp1
\advance\dimen0by5\baselineskip
\advance\dimen0by0.4true cm
\ifdim\dimen0>\vsize\pageinsert\box1\vfill\endinsert
\else
\if Y\FIG\stackfigbox\else
\dimen0=\pagetotal\ifdim\dimen0<\pagegoal
\advance\dimen0by\ht1\advance\dimen0by\dp1\advance\dimen0by1.4true cm
\ifdim\dimen0>\pagegoal\stackfigbox
\else\box1\vskip4true mm\fi
\else\box1\vskip4true mm\fi\fi\fi}
%
\def\begfig#1cm#2\endfig{\par
\setbox1=\vbox{\dimen0=#1true cm\advance\dimen0
by1true cm\kern\dimen0#2}\placefigure}
\def\begdoublefig#1cm #2 #3 \enddoublefig{\begfig#1cm%
\raggedright\hyphenpenalty=10000
\vskip-.8333\baselineskip\line{\vtop{\hsize=0.46\hsize#2}\hfill
\vtop{\hsize=0.46\hsize#3}}\endfig}
\def\figure#1#2{\vskip1true cm
\setbox0=\vbox{\noindent\petit{\bf
Fig.\ts#1\unskip.\ }\ignorespaces #2\smallskip
\count255=0\global\advance\count255by\prevgraf}%
\ifnum\count255>1\box0\else
\leftline{\petit{\bf Fig.\ts#1\unskip.\
}\ignorespaces#2}\smallskip\fi}

\def\begtab#1cm#2\endtab{\par
   \ifvoid\topins\midinsert\medskip\vbox{#2\kern#1true cm}\endinsert
   \else\topinsert\vbox{#2\kern#1true cm}\endinsert\fi}
\def\begpet{\vskip6pt\bgroup\petit}
\def\endpet{\vskip6pt\egroup}
\newdimen\refindent
\def\begref#1{\titlea{}{\ignorespaces#1}%
\bgroup\petit
\setbox0=\hbox{99.\enspace}\refindent=\wd0\relax}
\def\refno#1{\goodbreak
\hangindent\refindent\hangafter=1\noindent\hbox
to\refindent{\hss\ignorespaces#1\unskip\enspace}\ignorespaces}
\def\endref{\goodbreak\endpet}
\def\vec#1{{\textfont1=\tenbf\scriptfont1=\sevenbf
\textfont0=\tenbf\scriptfont0=\sevenbf
\mathchoice{\hbox{$\displaystyle#1$}}{\hbox{$\textstyle#1$}}
{\hbox{$\scriptstyle#1$}}{\hbox{$\scriptscriptstyle#1$}}}}
\def\petit{\def\rm{\fam0\ninerm}%
\textfont0=\ninerm \scriptfont0=\sixrm \scriptscriptfont0=\fiverm
 \textfont1=\ninei \scriptfont1=\sixi \scriptscriptfont1=\fivei
 \textfont2=\ninesy \scriptfont2=\sixsy \scriptscriptfont2=\fivesy
 \def\it{\fam\itfam\nineit}%
 \textfont\itfam=\nineit
 \def\sl{\fam\slfam\ninesl}%
 \textfont\slfam=\ninesl
 \def\bf{\fam\bffam\ninebf}%
 \textfont\bffam=\ninebf \scriptfont\bffam=\sixbf
 \scriptscriptfont\bffam=\fivebf
 \def\sans{\fam\sansfam\ninesans}%
 \textfont\sansfam=\ninesans \scriptfont\sansfam=\sixsans
 \scriptscriptfont\sansfam=\fivesans
 \def\tt{\fam\ttfam\ninett}%
 \textfont\ttfam=\ninett
 \normalbaselineskip=10pt
 \setbox\strutbox=\hbox{\vrule height7pt depth2pt width0pt}%
 \normalbaselines\rm
\def\vec##1{{\textfont1=\ninebf\scriptfont1=\sixbf
\textfont0=\ninebf\scriptfont0=\sixbf
\mathchoice{\hbox{$\displaystyle##1$}}{\hbox{$\textstyle##1$}}
{\hbox{$\scriptstyle##1$}}{\hbox{$\scriptscriptstyle##1$}}}}}
\nopagenumbers
%
\let\header=Y
\let\FIG=N
\newbox\figbox
\output={\if N\header\headline={\hfil}\fi\plainoutput
\global\let\header=Y
\if Y\FIG\topinsert\unvbox\figbox\endinsert\global\let\FIG=N\fi}
\let\lasttitle=N
\def\centerpar#1{{\parfillskip=0pt
\rightskip=0pt plus 1fil
\leftskip=0pt plus 1fil
\advance\leftskip by\oldparindent
\advance\rightskip by\oldparindent
\def\newline{\break}%
\noindent\ignorespaces#1\par}}
\let\INS=N%
\headline={\petit\def\newline{ }\def\fonote#1{}\ifodd\pageno
\rightheadline\else\leftheadline\fi}
\def\rightheadline{\hfil\the\AUTHOR\hfil\llap{\folio}}%
\def\leftheadline{\rlap{\folio}\hfil\the\AUTHOR\hfil}%
\nopagenumbers
\let\header=Y
\def\makestars{\count255=0\relax$^{
\loop\star\advance\count255 by 1\relax
\ifnum\count255<\footcount\repeat}$}
\let\lasttitle=N
 \def\head#1{\ifodd\pageno\else\null\vfill\supereject\fi
 \null\vskip2.5cm
 \let\header=N\bgroup
 \textfont0=\tafontt \scriptfont0=\tafonts \scriptscriptfont0=\tafontss
 \textfont1=\tamt \scriptfont1=\tams \scriptscriptfont1=\tams
 \textfont2=\tast \scriptfont2=\tass \scriptscriptfont2=\tasss
 \par\baselineskip=16pt
     \lineskip=16pt
     \tafontt
\def\textindent##1{\noindent\hbox
to0.5\oldparindent{##1\hss}\ignorespaces}%
\def\fonote##1{\advftncnt\makestars\begingroup\petit
\parfillskip=0pt plus 1fil
\vfootnote{\makestars}{##1\vskip-9.69pt}\endgroup}
     \pretolerance=10000
     \rightskip=0pt plus 1fil
     \noindent
     \ignorespaces#1
     \vskip16pt
     \egroup
     \nobreak
     \parindent=0pt
     \everypar={\global\parindent=1.5em
     \global\let\lasttitle=N\global\everypar={}}%
     \global\let\lasttitle=H%
     \ignorespaces}
 \def\subhead#1{
 \bgroup
 \textfont0=\tbfontt \scriptfont0=\tbfonts \scriptscriptfont0=\tbfontss
 \textfont1=\tbmt \scriptfont1=\tbms \scriptscriptfont1=\tbmss
 \textfont2=\tbst \scriptfont2=\tbss \scriptscriptfont2=\tbsss
 \par\baselineskip=14pt
     \lineskip=14pt
     \tbfontt
\def\textindent##1{\noindent\hbox
to0.5\oldparindent{##1\hss}\ignorespaces}%
\def\fonote##1{\advftncnt\makestars\begingroup\petit
\parfillskip=0pt plus 1fil
\vfootnote{\makestars}{##1\vskip-9.69pt}\endgroup}
     \pretolerance=10000
     \rightskip=0pt plus 1fil
     \vskip-6pt
     \noindent
     \ignorespaces#1
     \vskip16pt\egroup
     \nobreak
     \parindent=0pt
     \everypar={\global\parindent=1.5em
     \global\let\lasttitle=N\global\everypar={}}%
     \global\let\lasttitle=S%
     \ignorespaces}
\newtoks\AUTHOR
\AUTHOR={Missing surname(s) of the author(s), use $\backslash${\tt
authorrunning}}
\def\authorrunning#1{\global\AUTHOR={\ignorespaces#1\unskip}}
\def\author#1{
\global\footcount=0
\bgroup
\pretolerance=10000
\noindent\ignorespaces#1\vskip2\baselineskip\egroup\let\lasttitle=X}
\def\address{\vskip2\baselineskip\bgroup\petit
\parskip=0pt
\parindent=0pt
\obeylines\addrreess}
\def\addrreess#1{\noindent\ignorespaces#1\bigskip\egroup}
\def\abstract#1{\noindent
\ignorespaces#1\vskip24pt
\let\lasttitle=X
\parindent=0pt
\everypar={\global\parindent=\oldparindent
\global\let\lasttitle=N\global\everypar={}}%
\ignorespaces}
 \def\titlea#1#2{\if N\lasttitle\else\vskip-24pt
     \fi
     \vskip24pt plus 4pt minus4pt
     \bgroup
\textfont0=\tenbf \scriptfont0=\sevenbf \scriptscriptfont0=\fivebf
\textfont1=\tams \scriptfont1=\tamss \scriptscriptfont1=\tbmss
     \lineskip=0pt
     \pretolerance=10000
     \noindent
     \bf
     \rightskip 0pt plus 6em
     \setbox0=\vbox{\vskip23pt\def\fonote##1{}%
     \noindent
     \if!#1!\ignorespaces#2
     \else\setbox0=\hbox{\ignorespaces#1\unskip.\enspace}%
     \hangindent=\wd0
     \hangafter=1\box0\ignorespaces#2\fi
     \vskip10pt}%
     \dimen0=\pagetotal\advance\dimen0 by-\pageshrink
     \ifdim\dimen0<\pagegoal
     \dimen0=\ht0\advance\dimen0 by\dp0\advance\dimen0 by
     3\normalbaselineskip
     \advance\dimen0 by\pagetotal
     \ifdim\dimen0>\pagegoal\eject\fi\fi
     \noindent
     \if!#1!\ignorespaces#2
     \else\setbox0=\hbox{\ignorespaces#1\unskip.\enspace}%
     \hangindent=\wd0
     \hangafter=1\box0\ignorespaces#2\fi
     \vskip12pt plus4pt minus4pt\egroup
     \nobreak
     \parindent=0pt
     \everypar={\global\parindent=\oldparindent
     \global\let\lasttitle=N\global\everypar={}}%
     \global\let\lasttitle=A%
     \ignorespaces}
 \def\titleb#1#2{\if N\lasttitle\else\vskip-24pt
     \fi
     \vskip24pt plus 4pt minus4pt
     \bgroup
     \it
     \lineskip=0pt
     \pretolerance=10000
     \noindent
     \rightskip 0pt plus 6em
     \setbox0=\vbox{\vskip23pt\def\fonote##1{}%
     \noindent
     \if!#1!\ignorespaces#2
     \else\setbox0=\hbox{\ignorespaces#1\unskip.\enspace}%
     \box0%
     \ignorespaces#2\fi
     \vskip6pt}%
     \dimen0=\pagetotal\advance\dimen0 by-\pageshrink
     \ifdim\dimen0<\pagegoal
     \dimen0=\ht0\advance\dimen0 by\dp0\advance\dimen0 by
     2\normalbaselineskip
     \advance\dimen0 by\pagetotal
     \ifdim\dimen0>\pagegoal\eject\fi\fi
     \noindent
     \if!#1!\ignorespaces#2
     \else\setbox0=\hbox{\ignorespaces#1\unskip.\enspace}%
     \box0%
     \ignorespaces#2\fi
     \vskip12pt plus4pt minus4pt\egroup
     \nobreak
     \parindent=0pt
     \everypar={\global\parindent=\oldparindent
     \global\let\lasttitle=N\global\everypar={}}%
     \global\let\lasttitle=B%
     \ignorespaces}
 \def\titlec#1{\if N\lasttitle\else\vskip-\baselineskip
     \fi
     \vskip12pt plus 4pt minus 4pt
     \bgroup
     \it
     \noindent
     \ignorespaces#1\unskip.\ \egroup
     \ignorespaces}
 \def\titled#1{\if N\lasttitle\removelastskip\vskip\baselineskip
     \fi
     \bgroup
     \noindent
\textfont0=\tenbf \scriptfont0=\sevenbf \scriptscriptfont0=\fivebf
\textfont1=\tams \scriptfont1=\tamss \scriptscriptfont1=\tbmss
     \tenbf
     \ignorespaces#1\unskip.\ \egroup
     \ignorespaces}
\let\ts=\thinspace
\def\footnoterule{\kern-3pt\hrule width 2true cm\kern2.6pt}
\newcount\footcount \footcount=0
\def\advftncnt{\advance\footcount by1\global\footcount=\footcount}
\def\fonote#1{\advftncnt$^{}$\begingroup\petit
\parfillskip=0pt plus 1fil
\def\textindent##1{\noindent\hbox
to\oldparindent{##1\hss}\ignorespaces}%
\vfootnote{$^{}$}{#1\vskip-9.69pt}\endgroup}
\newdimen\itemindent
\setbox0=\hbox{99.\enspace}\itemindent=\wd0\relax
\def\item#1{\par\hangindent\itemindent\noindent\hbox
to\itemindent{\hss\ignorespaces#1\unskip\enspace}%
\ignorespaces}
\newdimen\itemitemindent
\itemitemindent=\itemindent
\advance\itemitemindent by\oldparindent\relax
\def\itemitem#1{\par\hangindent\itemitemindent\noindent
\kern\itemindent\hbox
to\oldparindent{\hss\ignorespaces#1\unskip\enspace}%
\ignorespaces}
\def\appendix#1{\titlea{}{Appendix #1}}
\def\acknow#1{
 \if N\lasttitle\removelastskip\vskip\baselineskip
     \fi
     \bgroup\petit
     \noindent
     {\it Acknowledgement.\ }%
     \ignorespaces#1\par\egroup
     \ignorespaces}
\def\newenvironment#1#2#3#4{\long\def#1##1##2{\removelastskip
\vskip\baselineskip\noindent{#3#2\if!##1!.\else\ ##1.\fi
\ }{#4\ignorespaces##2}\vskip\baselineskip}}
\newenvironment\lemma{Lemma}{\tenbf}{\it}
\newenvironment\proposition{Proposition}{\tenbf}{\it}
\newenvironment\theorem{Theorem}{\tenbf}{\it}
\newenvironment\corollary{Corollary}{\tenbf}{\it}
\newenvironment\example{Example}{\tenbf}{\rm}
\newenvironment\exercise{Exercise}{\tenbf}{\rm}
\newenvironment\problem{Problem}{\tenbf}{\rm}
\newenvironment\solution{Solution}{\tenbf}{\rm}
\newenvironment\definition{Definition}{\tenbf}{\rm}
\newenvironment\remark{Remark}{\tenbf}{\rm}
\newenvironment\note{Note}{\tenbf}{\rm}
\long\def\proof{\removelastskip\vskip\baselineskip\noindent{\it
Proof.\quad}\ignorespaces}
\def\frame#1{\bigskip\vbox{\hrule\hbox{\vrule\kern5pt
\vbox{\kern5pt\advance\hsize by-10.8pt
\centerline{\vbox{#1}}\kern5pt}\kern5pt\vrule}\hrule}\bigskip}
\def\frameddisplay#1#2{$$\vcenter{\hrule\hbox{\vrule\kern5pt
\vbox{\kern5pt\hbox{$\displaystyle#1$}%
\kern5pt}\kern5pt\vrule}\hrule}\eqno#2$$}
\def\typeset{\petit\noindent This article was processed by the authors
\newline using the Springer-Verlag \TeX{} mamath macro package
1990.\par}
\outer\def\bye{\bigskip\bigskip\typeset
\footcount=1\ifx\speciali\undefined\else
\loop\smallskip\noindent special character No\number\footcount:
\csname special\romannumeral\footcount\endcsname
\advance\footcount by 1\global\footcount=\footcount
\ifnum\footcount<11\repeat\fi
\vfill\supereject\end}

\catcode`\@=11
\font\smc=cmcsc9
\newfam\ssffam
\font\tenssf=cmss10
\font\eightssf=cmss8
\def\ssf{\fam\ssffam}
  \textfont\ssffam=\tenssf \scriptfont\ssffam=\eightssf
\def\sf#1{\leavevmode\skip@\lastskip\unskip\/%
       \ifdim\skip@=\z@\else\hskip\skip@\fi{\ssf#1}}
\def\rom#1{\leavevmode\skip@\lastskip\unskip\/%
        \ifdim\skip@=\z@\else\hskip\skip@\fi{\rm#1}}
\def\newmcodes@{\mathcode`\'"27\mathcode`\*"2A\mathcode`\."613A%
 \mathcode`\-"2D\mathcode`\/"2F\mathcode`\:"603A }
\def\operatorname#1{\mathop{\newmcodes@\kern\z@\fam\z@#1}}
\catcode`\@=\active

\def\N{{\sf N}}

\def\gp#1{\langle#1\rangle}
\def\m1{^{-1}}
\def\ov1{\overline}
\def\supp{\operatorname{supp\,}}
\def\,{\thinspace}
\let\tsize\textstyle
\def\text#1{\hbox{\rm#1}}

\relpenalty=10000
\binoppenalty=10000
\mathsurround=1pt

\def\cite#1{\rom{[#1]}}
\def\aab{1}
\def\bgs{2}
\def\bsb{3}
\def\bsn{4}
\def\vab{5}
\def\mhs{6}
\def\hig{7}
\def\kap{8}
\def\mfn{9}
\def\nob{10}
\def\rse{11}
\hyphenation{manu-scripta}

\input amssym.def
\input amssym.tex
\magnification 1000 
\head{Unitary units in modular group algebras}
\author{Victor Bovdi and L.\,G. Kov\'acs
\fonote
{Research partly supported by the Hungarian National Foundation for Scientific
Research grant no. T4265.\endgraf
 The second author is indebted to the `Universitas' Foundation and the Lajos
Kossuth University of Debrecen, Hungary, for warm hospitality and generous
support during the period when this work began.}}
\authorrunning{Bovdi and Kov\'acs}

\abstract
{Let $p$ be a prime, $K$ a field of characteristic $p$, $G$ a locally finite
$p$-group, $KG$ the group algebra, and $V$ the group of the units of $KG$ with
augmentation $1$. The anti-automorphism $g\mapsto g\m1$ of $G$ extends linearly
to $KG$; this extension leaves $V$ setwise invariant, and its restriction to $V$
followed by $v\mapsto v\m1$ gives an automorphism of $V$. The elements of $V$
fixed by this automorphism are called {\it unitary\/}; they form a subgroup.
Our first theorem describes the $K$ and $G$ for which this subgroup is normal
in $V$.
\endgraf
 For each element $g$ in $G$, let $\ov1g$ denote the sum (in $KG$) of the
distinct powers of $g$. The elements $1+(g-1)h\ov1g$ with $g,\,h\in G$ are the
{\it bicyclic} units of $KG$. Our second theorem describes the $K$ and $G$ for
which all bicyclic units are unitary.}

\titlea{1}{Introduction}
 Let $KG$ be the group algebra of a group $G$ over a commutative ring $K$
(with~$1$) and $V(KG)$ the group of normalized units (that is, of the units with
augmentation $1$) in $KG$. The anti-automorphism $g\mapsto g\m1$ extends
linearly to an anti-automorphism $a\mapsto a^*$ of $KG$; this extension leaves
$V(KG)$ setwise invariant, and its restriction to $V(KG)$ followed by
$v\mapsto v\m1$ gives an automorphism of $V(KG)$. The elements of $V(KG)$ fixed
by this automorphism are the {\it unitary normalized units} of $KG$; they form a
subgroup which we denote by $V_*(KG)$. (Interest in unitary units arose in
algebraic topology, and a more general definition, involving an `orientation
homomorphism', is also current; the special case we use here arises when the
orientation homomorphism is trivial.)

 The first question considered here is to find the pairs $K$, $G$ for which
$V_*(KG)$ is normal in $V(KG)$. (Since each unit of a group algebra is a scalar
multiple of a normalized unit, if $V_*(KG)$ is normal in $V(KG)$ then it is
normal also in the group of all units of $KG$.) For $K=\Bbb Z$, this question
was discussed (without any restriction on the orientation homomorphism) by
A.\,A. Bovdi and S.\,K. Sehgal in \cite{\bsn}. Here we deal with the `modular'
case, that is, with the case of $K$ a field of prime characteristic $p$ and $G$
a locally finite $p$-group.

\theorem{1.1}
 {Let $K$ be a field of prime characteristic $p$ and let $G$ be a nonabelian
locally finite $p$-group. The subgroup $V_*(KG)$ is normal in $V(KG)$ if and
only if $p=2$ and $G$ is the direct product of an elementary abelian group with
a group $H$ for which one of the following holds\rom:
\itemitem{\rm(i)} $H$ has no direct factor of order $2$, but it is a semidirect
product of a group $\gp{h}$ of order $2$ and an abelian $2$-group $A$, with
$h\m1ah=a\m1$ for all $a$ in $A$\rom;
\itemitem{\rm(ii)} $H$ is an extraspecial $2$-group, or the central product of
such a group with a cyclic group of order $4$.}

 We work with the definition that a $p$-group is extraspecial if its centre,
commutator subgroup and Frattini subgroup are equal and have order $p$: we do
not require the group itself to be finite.

 The proof of this theorem will be given in Section 2. The reason we take the
$p$-group $G$ locally finite is that, as is well known, this ensures that each
non-unit of $KG$  lies in the augmentation ideal.

 Every group $G$ may be written (see Lemma 2.3) as a direct product of an
elementary abelian $2$-group $E$ and a group $H$ which has no direct factor of
order $2$ (we do not exclude $E=1$ or $H=1$). The isomorphism type of $G$
determines the isomorphism types of $E$ and $H$, and vice versa.

 It is easy to verify that if $H$ satisfies (i) then
$A=\gp{\,a\in H\mid a^2\neq1\,}$ and $A$ has no direct factor of order $2$.
Conversely, if $A$ is a nontrivial abelian $2$-group without a direct factor of
order $2$ and $H$ is formed as the semidirect product indicated, then (i) holds.
The classification of the groups $H$ of this kind is thus reduced to the
classification of abelian $2$-groups, a problem whose solution in terms of Ulm
invariants is well known in the finite or countably infinite case but is beyond
reach in general.

 As to case (ii), the classification of finite extraspecial groups is well
known. Equally conclusive results were obtained for extraspecial groups of
countably infinite order by M.\,F. Newman in \cite{\mfn}; he also showed
there that no such results can be expected for extraspecial groups of arbitrary
order.

 The only group $H$ which satisfies both conditions (i) and (ii) is the dihedral
group of order $8$.

\smallskip

 The second part of the paper concerns the bicyclic units introduced in Ritter
and Sehgal \cite{\rse}. For $K$ a commutative ring and $g$ an element of finite
order $|g|$ in a group $G$, let $\ov1g$ denote the sum (in $KG$) of the distinct
powers of $g$:
$$
\ov1g=\sum_{i=0}^{|g|-1} g^i.
$$
 If also $h\in G$, put
$$
 u_{g,h}=1+(g-1)h\ov1g.
$$
 Note that $1-(g-1)h\ov1g$ is a two-sided inverse for $u_{g,h}$ and the
augmentation of $u_{g,h}$ is $1$, so $u_{g,h}$ is a normalized unit. The
elements of this form are called {\it bicyclic units.}

 The problem considered here is to find the $K$ and $G$ for which each bicyclic
unit of $KG$ is unitary. It is easy to see that $u_{g,h}=1$ if and only if the
cyclic group $\gp{g}$ is normalized by $h$, and that if $K=\Bbb Z$ then $1$ is
the only bicyclic unit which is unitary. Thus the $G$ which can partner
$K=\Bbb Z$ are precisely the groups in which every subgroup of finite order is
normal. (The situation is not so simple for $K=\Bbb Z$ when unitarity is defined
with reference to a nontrivial orientation homomorphism: see A.\,A. Bovdi and
S.\,K. Sehgal \cite{\bsb}.) In private communication, Professor Sehgal has
directed attention to the modular case: what can one say when $K$ is of
characteristic $p$ and $G$ is a $p$-group?

\theorem{1.2}
 {Let $p$ be a prime, $K$ a commutative ring of characteristic $p$, and $G$ a
nonabelian $p$-group. All bicyclic units of $KG$ are unitary if and only if
$p=2$ and $G$  is the direct product of an elementary abelian group and a
group $H$ for which one of the following holds\rom:
\itemitem{\rm(i)} $H$ has an abelian subgroup $A$ of index $2$ such that
conjugation by an element of $H$ outside $A$ inverts each element of $A$\rom;
\itemitem{\rm(ii)} $H$ is an extraspecial $2$-group, or the central product of
such a group with a cyclic group of order $4$\rom;
\itemitem{\rm(iii)} $H$ is the direct product of a quaternion group of order $8$
and a cyclic group of order $4$, or the direct product of two quaternion groups
of order $8$\rom;
\itemitem{\rm(iv)} $H$ is the central product of the group
$\gp{\,x,y\mid x^4=y^4=1,\ x^2=[y,x]\,}$  with a quaternion group of order $8$,
the nontrivial element common to the two central factors being $x^2y^2$\rom;
\itemitem{\rm(v)} $H$ is isomorphic to one of the groups $H_{32}$ and $H_{245}$
defined below.}

\noindent
 The relevant definitions are:
$$
\eqalign{
 H_{32}=\bigl\langle\,x,y,u\bigm|\ &x^4=y^4=1,\cr  &x^2=[y,x],\cr
&y^2=u^2=[u,x],\cr  &x^2y^2=[u,y]\,\bigr\rangle,\cr
\noalign{\vskip2pt}\allowbreak
 H_{245}=\bigl\langle\,x,y,u,v\bigm|\ &x^4=y^4=[v,u]=1,\cr
&x^2=v^2=[y,x]=[v,y],\cr  &y^2=u^2=[u,x],\cr
&x^2y^2=[u,y]=[v,x]\,\bigr\rangle.\cr}
$$
 The subscripts are the serial numbers of these groups in the CAYLEY library of
groups of order dividing $128$ described by Newman and O'Brien in \cite{\nob}.
It is a mere coincidence that $H_{32}$ has order $32$. The other group,
$H_{245}$, is one of the two Suzuki $2$-groups (see Higman \cite{\hig}) of order
$64$. That CAYLEY library provides us not only with serial numbers but also
with the final step in the proof of Theorem 1.2: we are indebted to Dr E.\,A.
O'Brien for extracting the list (Lemma 4.1) of the groups $H$ of order
dividing $128$ whose Frattini subgroup is central, noncyclic, of order $4$, and
contains all elements of order $2$ in $H$. His list shows that all groups of
this kind have order dividing $64$, so they may also be found in Hall and Senior
\cite{\mhs}, where $H_{32}$ and $H_{245}$ are labelled $32\,\Gamma_4c_3$ and
$64\,\Gamma_{13}a_5$. (Of course, without the much larger CAYLEY library we
could not see that no group of order $128$ satisfies these criteria.)

 Several of the comments we made after Theorem 1.1 apply here as well. In (i)
we can assume that $H$ has no direct factor of order $2$: the isomorphism type
of $H$ is then determined by $G$. We may also assume there that $|H|>8$, for the
two nonabelian groups of order $8$ occur also under (ii). It is easy to see that
then $H$ has only one abelian  subgroup of index $2$, so the isomorphism type of
$A$ is in turn determined by $H$. Moreover, the squares of all the elements of
$H$ outside $A$ are equal to each other, and this element, $a_0$ say, has
order at most $2$. Of course the height of $a_0$ in $A$ is also an isomorphism
invariant of $H$. Mackey's proof of Ulm's Theorem (given in Kaplansky
\cite{\kap}) shows that if two countable $p$-groups have the same Ulm invariants
and we are given a height-preserving isomorphism from a finite subgroup of one
to a subgroup of the other, then this will extend to an isomorphism of the two
groups. It follows that in the finite or countably infinite case the Ulm
invariants of $A$ together with the height in $A$ of the common square of the
elements outside $A$ form a {\it complete} set of invariants for $H$.
 Conversely,
if $a_0$ is any element of order at most $2$ in an abelian $2$-group $A$ with
$|A|>4$, then the group $H$ defined by
$$
 H=\bigl\langle\,A,\,h\bigm| h^2=a_0,\ h\m1ah=a\m1 \text{ for all }\,a\,
\text{ in }\,A\,\bigr\rangle
$$
 satisfies (i) and $|H|>8$. (The reader may like to work out how the relevant
invariants must be restricted to ensure that $H$ is nonabelian and has no direct
factor of order $2$.)

 The proof of Theorem 1.2 splits naturally into a ring-theoretic part and a
group-theoretic part, which are presented in Section 3 and Section 4. Their
conclusions  make sense separately and may be of some independent interest, so
we state them here.

\lemma{1.3}
 {Let $K$ be a commutative ring with $1$ and $G$ a group. Suppose that $g$ is an
element of finite order in $G$ and $h$ is an element of $G$ which does not
normalize $\gp{g}$. The bicyclic unit $u_{g,h}$ is unitary if and only if the
characteristic of $K$ is $2$ while $\gp{g^2}$ is normalized by $h$ and contains
either $h^2$ or $(hg)^2$.}

 (In particular, this confirms the comment above that in $\Bbb ZG$ the only
unitary bicyclic unit is $1$. In fact, $V_*(\Bbb ZG)$ is always $G$ itself, as
was observed in Lemma~1 of A.\,A. Bovdi \cite{\aab}.)

\lemma{1.4}
 {Let $G$ be a nonabelian $2$-group such that if $g,h\in G$ then $\gp{g^2}$ is
normal in $G$ and  $\gp{g,h}/\gp{g^2}$ is either abelian or dihedral. Then $G$
is the direct product of an elementary abelian group with a nonabelian group $H$
for which one of the conditions {\rm (i)--(v)} of Theorem \rom{1.2} holds.
Conversely, every $2$-group $G$ of this kind satisfies our hypotheses.}

\titlea{2}{Normal unitary subgroup}
 The aim of this section is to prove Theorem 1.1. It will be convenient to use
$y^*=y\m1$ as the test of whether a unit $y$ is unitary. At first, $K$ can be
any commutative ring with $1$ and $G$ any group.

\lemma{2.1}
 {For $x\in V(KG)$ and $y\in V_*(KG)$, we have $x\m1yx\in V_*(KG)$ if and only
if $xx^*$ commutes with $y$.}

\proof
 Clearly, $(x\m1yx)^*=(x\m1yx)\m1$ means that $x^*y^*(x^*)\m1=x\m1y\m1x$ which
in turn is equivalent to $xx^*y^*=y\m1xx^*$. As we are given that $y^*=y\m1$,
this proves the lemma.\qed

\smallskip

 As $G\leq V_*(KG)$, an element which commutes with every element of $V_*(KG)$
is central in $KG$. Thus Lemma 2.1 gives the following.

\corollary{2.2}
 {The subgroup $V_*(KG)$ is normal in $V(KG)$ if and only if all elements
of the form $xx^*$ with $x\in V(KG)$ are central in $KG$. \qed}

 (For the case of $K=\Bbb Z$, this is a special case of Lemma 2 of \cite{\bsn};
the proof we have given comes from that paper.)

\medskip

\noindent{\it Proof of Theorem \rom{1.1}.\quad}From the simple fact that over a
field of characteristic $p$ a finite $p$-group has only one irreducible
representation, it follows readily that under the hypotheses of the theorem the
augmentation ideal of $KG$ is locally nilpotent and so each element outside that
ideal is a unit. Differently put, if $x$ is a non-unit in $KG$ then $1+KG\in
V(KG)$.

 Suppose first that $V_*(KG)$ is normal in $V(KG)$. By Corollary 1, if $x$ is
a normalized unit then $xx^*$ is central. As each unit is a scalar multiple
of a normalized unit, the same conclusion is available whenever $x$ is a
unit. It follows that if $x$ is any unit then $xx^*=x\m1(xx^*)x=x^*x$. If $x$
is a non-unit in $KG$, then $1+x$ is a unit and so $(1+x)(1+x)^*=(1+x)^*(1+x)$,
whence again $xx^*=x^*x$.

 A group algebra in which $xx^*=x^*x$ holds for every element $x$ is called
{\it normal.}  A.\,A.\ Bovdi, P.\,M.\ Gudivok and M.\,S.\ Semirot proved in
\cite{\bgs} that the group algebra of a nonabelian group $G$ over a commutative
ring $K$ is normal if and only if either $G$ is hamiltonian or the
characteristic of $K$ is $2$ and $G$ is a direct product of an elementary
abelian $2$-group with a group $H$ such that (i) or (ii) holds. Thus the
proof of our `only if' claim is complete.

 Suppose next that $p=2$ and $G=E\times H$ with $E$ elementary abelian and $H$
satisfying (i) or (ii). In view Corollary 2.2, what we have to show is that
$xx^*$ is central whenever $x\in V(KG)$.

 Consider first the case (i). Then each element $x$ of $KG$ can be written as
$x=x_1+x_2h$ with $x_1,\,x_2\in K(E\times A)$, and of course
$hx_ih=x_i^*$. Using again that $h^2=1$, that $K(E\times A)$ is commutative, and
that the characteristic is $2$, we see that
$xx^*=(x_1+x_2h)(x_1^*+hx_2^*)=x_1x_1^*+2x_1x_2h+x_2x_2^*=x_1x_1^*+x_2x_2^*$.
Thus $xx^*$ lies in the commutative algebra $K(E\times H)$ and is easily seen to
commute with $h$, so it is central in $KG$.

 Consider next the case (ii). Then the commutator subgroup of $G$ has
only one nontrivial element; call that $c$, and write $I$ for the ideal of $KG$
generated by $1+c$. This element $c$ is central in $G$, while if $g,\,h\in G$
then either $hg=gh$ or $hg=ghc$: so either $hg(1+c)=gh(1+c)=g(1+c)h$ or
$hg(1+c)=ghc(1+c)=gh(1+c)=g(1+c)h$ proves that $g(1+c)$ commutes with $h$. It
follows that every element of $I$, and therefore also every element of $1+I$, is
central in $KG$.

 Let $\gamma$ be the natural homomorphism of $KG$ onto $K(G/\gp{c})$ defined
by $g\gamma=g\gp{c}$ for all $g$ in $G$. Of course $\gamma$ intertwines the
augmentation maps of the two group algebras, so if $x$ is a normalized unit in
$KG$ then $x\gamma$ is a normalized unit in $K(G/\gp{c})$. Further, $\gamma$
intertwines the anti-automorphism $*$ of $KG$ with the similarly defined
anti-automorphism of $K(G/\gp{c})$; we shall use $*$ also for the latter
anti-automorphism. Note that $K(G/\gp{c})$ is elementary abelian. It is an easy
exercise to see that in a characteristic $2$ group algebra of an elementary
abelian $2$-group each normalized unit is unitary. In particular, if
$x\in V(KG)$ then $x\gamma$ is unitary, so
$(xx^*)\gamma=(x\gamma)(x^*\gamma)=(x\gamma)(x\gamma)^*=1$, that is,
$xx^*\in1+\ker\gamma$. Since $I$ is minimal among the ideals for which
$c\equiv 1\,\bmod I$, it is precisely $\ker\gamma$. We have proved that
$xx^*\in 1+I$. By the conclusion of the previous paragraph, $xx^*$ is therefore
central in $KG$. The proof of the theorem is now complete. \qed

\smallskip

\noindent{\it Remarks.} \
 On any group algebra of an elementary abelian $2$-group, the
anti-automorphism $*$ is in fact just the identity map.

 For a generalization of the result of \cite{\bgs}, see \cite{\vab}.

\medskip

 We conclude this section with a purely group-theoretic lemma which was
mentioned in the introduction's comments on Theorem 1.1.

\lemma{2.3}
  {Every group $G$ is a direct product $E\times H$ of an elementary abelian
$2$-group $E$ and a group $H$ which has no direct factor of order $2$. If
$G=E_1\times H_1$ is another such decomposition of $G$, then $E_1\cong E$ and
$H_1\cong H$.}

\proof
  Let $\sf Z(G)$ denote the centre of $G$; set
$A=\sf A(G)=\gp{\,a\in \sf Z(G)\mid a^2=1\,}$ and
$B=\sf B(G)=\gp{\,g^2\mid g\in G\,}$. Let $E$ be a direct complement to
$A\cap B$ in $A$, and $H/B$ a direct complement to $AB/B$ in $G/B$: then
$$
 G=ABH=EBH=EH \quad\text{while}\quad E\cap H\leq E\cap AB\cap H=E\cap B=1,
$$
 so $G=E\times H$. Here $E$ is elementary abelian because $A$ is. If
$H=C\times K$ with $|C|\leq2$, then $G=E\times C\times K$ and $B=\sf B(K)$, so
$C\leq A=E\times(A\cap B)\leq E\times K$ yields that $C=1$. This proves that
$H$ has no direct factor of order $2$, that is, $\sf A(H)\leq\sf B(H)$. If
$G=E_1\times H_1$ is another decomposition with $E_1$ elementary abelian and
$\sf A(H_1)\leq\sf B(H_1)$, then $E_1\cong A/(A\cap B)\cong E$, and $H_1/B$ is
another direct complement to $AB/B$ in $G/B$: so we also have $G=E\times H_1$
and therefore $H_1\cong G/E\cong H$.\qed

\titlea{3}{Unitary bicyclic units}
The aim of this section is to prove Lemma 1.3. Accordingly, $K$ is once again
an arbitrary commutative ring with 1 and $G$ is an arbitrary group. Recall the
definition $u_{g,h}=1+(g-1)h\ov1g$, and note that $u_{g,h}=u_{g,hg}\in V(KG)$,
with
$$
 u_{g,h}\m1=1-(g-1)h\ov1g\text{\qquad and\qquad}u_{g,h}^*=1+\ov1gh\m1(g\m1-1).
$$

 The {\it support} of an element $a$ of $KG$ is the set of those elements of
$G$ which occur with nonzero coefficient in the expression of $a$ as $K$-linear
combination of elements of $G$:
$$
\tsize
\supp\sum_{g\in G}\alpha_gg=\{\,g\in G\mid\alpha_g\neq0\,\}.
$$
 Two simple observations about bicyclic units will be used without reference.
First, if $h$ normalizes $\gp{g}$ then $h\ov1g=\ov1gh$ and so $u_{g,h}=1$.
Second, if $h$ does not normalize $\gp{g}$ then $u_{g,h}\neq1$; indeed, in this
case no element of $G$ can occur more than once in the expansion of
$1+(g-1)h\ov1g$, so the support of $u_{g,h}$ has cardinality $1+|g|+|g|$.
Explicitly, if $h\notin\N(\gp{g})$ then
$$
\supp u_{g,h}=\{1\}\cup\bigl\{ghg^i\bigm|0\leq i<|g|\,\bigr\}\cup
\bigl\{\,hg^i\bigm|0\leq i<|g|\,\bigr\}.
$$
{\it Proof of Lemma \rom{1.3}.\quad}Suppose that $u_{g,h}^*=u_{g,h}\m1$.
If the characteristic of $K$ were not $2$, we could argue that in the expression
of $u_{g,h}\m1$ as $1-(g-1)h\ov1 g$ both $h$ and $hg$ have coefficient $1$ while
in that of $u_{g,h}^*$ the only nontrivial elements of $G$ with coefficient $1$
are the $g^ih\m1g\m1$, hence there exist $i$, $j$ such that $h=g^ih\m1g\m1$ and
$hg=g^jh\m1g\m1$, and then
$$
 hgh\m1=\bigl(g^jh\m1g\m1\bigr)\bigl(g^ih\m1g\m1\bigr)\m1
=g^{j-i}\in\gp{g}
$$
 contradicts the assumption that $h\notin\N(\gp{g})$. Thus the characteristic of
$K$ is $2$.

 Note that $|g^2|$ can only be $|g|$ or $|g|/2$. We exploit this repeatedly, for
it yields that once we show $h\m1g^2h\in\gp{g}$, it follows that $h$ normalizes
$\gp{g^2}$. Namely, if $h\m1g^2h\in\gp{g}$ then $h\notin\N(\gp{g})$ implies that
$\gp{h\m1g^2h}<\gp{g}$, whence $|h\m1g^2h|=|g^2|=|g|/2$, and then
$\gp{h\m1g^2h}$ must be the unique subgroup, $\gp{g^2}$, of index $2$ in
$\gp{g}$.

 Since $u_{g,h}\neq1$, the support of $u_{g,h}\m1$ is given by
$$
\supp u_{g,h}\m1=\{1\}\cup\bigl\{\,ghg^i\bigm|0\leq i<|g|\,\bigr\}\cup
\bigl\{\,hg^i\bigm|0\leq i<|g|\,\bigr\},
$$
while
$$
\supp u_{g,h}^*=\{1\}\cup\bigl\{\,g^ih\m1g\m1\bigm|0\leq i<|g|\,\bigr\}
\cup\bigl\{\,g^ih\m1\bigm|0\leq i<|g|\,\bigr\}.
$$
 Given our assumption that $u_{g,h}^*=u_{g,h}\m1$, these two supports are equal.
We now distinguish a number of cases according to the form in which various
elements of \,$\supp u_{g,h}\m1$ appear in \,$\supp u_{g,h}^*$.

 Suppose first that $h=g^ih\m1$, so $ghg=g^{i+1}h\m1g$. If $ghg=g^jh\m1$,
then  $h\m1gh=g^{j-i-1}\in\gp{g}$, contrary to $h\notin\N(\gp{g})$. Thus
$ghg=g^jh\m1g\m1$, and then $h\m1g^2h=g^{j-i-1}\in\gp{g}$, so $h$ normalizes
$\gp{g^2}$. Of course now $h^2=g^i\in\gp{g}$, so $h^2\notin\gp{g^2}$ would imply
that $\gp{g}=\gp{\gp{g^2},h^2}$, which is impossible because $h\in\N(\gp{g^2})$
but $h\notin\N(\gp{g})$. This proves that in this case $h^2\in\gp{g^2}$.

 Suppose next that $h=g^ih\m1g\m1$, and note that $h\notin\N(\gp{g})$
implies that $i\neq0$. If $|g|=2$ then this forces $i=1$, so conjugation by $g$
inverts $h$ and therefore $(hg)^2=1$. Suppose that $|g|>2$; then
$g^ih\m1g=hg^2\in\supp u_{g,h}\m1$. If $hg^2=g^jh\m1$, then
$h\m1gh=g^{j-i}\in\gp{g}$, contrary to $h\notin\N(\gp{g})$. Thus
$hg^2=g^jh\m1g\m1$, and then $h\m1g^2h=g^{j-i}\in\gp{g}$, so $h$ normalizes
$\gp{g^2}$. If $i$ is even, then in the factor group $\gp{g,h}/\gp{g^2}$ the
 image
of $g$ is the square of the image of $h$, but this is impossible because
$h\notin\N(\gp{g})$. So $i$ is odd, and then in $\gp{g,h}/\gp{g^2}$ conjugation
 by
the image of $g$ inverts the image of $h$ and therefore the image of $hg$
has order $2$: thus again $(hg)^2\in\gp{g^2}$.

 This completes the proof of the `only if' claim. For the proof of the `if'
claim, assume first that the characteristic of $K$ is $2$, and note that then
each bicyclic unit of $KG$ is its own inverse. Next, assume that $\gp{g^2}$ is
normalized by $h$ and contains either $h^2$ or $(hg)^2$. Since in any case
$u_{g,h}=u_{g,hg}$, we may assume without loss of generality that in fact
$h^2\in\gp{g^2}$. Since $\gp{g^2}$ is normalized by $h$ but $\gp{g}$ is not,
$|g|$ must be even, whence
$$
\ov1 g=(g+1)\ov1{g^2}=\ov1{g^2}(g\m1+1).
$$
 Further, as both $g$ and $h$ normalize $\gp{g^2}$, both commute with
$\ov1{g^2}$, while $h^2\in\gp{g^2}$ implies that $h^2\ov1{g^2}=\ov1{g^2}$ and so
$$
 h\ov1{g^2}=\ov1{g^2}h\m1.
$$
 Using again that the characteristic of $K$ is $2$, we can therefore argue that
$$
\eqalign{u_{g,h}\m1=u_{g,h}&=1+(g+1)h\ov1{g}\cr &=1+(g+1)h\ov1{g^2}(g\m1+1)\cr
&=1+(g+1)\ov1{g^2}h\m1(g\m1+1)\cr &=1+\ov1{g}h\m1(g\m1+1)\cr &=u_{g,h}^*\ ,}
$$
 as required.\qed

\titlea{4}{A certain class of groups}
 The rest of the paper will be taken up by the proof of Lemma 1.4.

 As usual, the Frattini subgroup of a group $H$ will be written $\Phi(H)$.
Recall that if $H$ is a finite $2$-group then
$\Phi(H)=\gp{\,h^2\mid h\in H\,}$. If $H$ is any $2$-group, we write
$\Omega(H)=\gp{\,h\in H\bigm|h^2=1\,}$.  The proof of Lemma 1.4 depends on a
result obtained for us by Dr~E.\,A.~O'Brien by inspecting the CAYLEY library
described in Newman and O'Brien \cite{\nob}.

\lemma{4.1 {\rm (O'Brien)}}
 {The groups $H$ of order dividing $128$ in which $\Phi(H)$ and $\Omega(H)$ are
equal, central, and of order $4$, are precisely the following: $C_4\times C_4$,
$C_4\rtimes C_4$, $C_4\rtimes Q_8$, and the groups named in parts
{\rm(iii)--(v)} of Theorem \rom{1.2}. \qed}

 Here $C_4$ and $Q_8$ stand for a cyclic group of order $4$ and a quaternion
group of order $8$,  while $C_4\rtimes C_4$ and $C_4\rtimes Q_8$ indicate
semidirect products which are not direct products (the last-named semidirect
factor not being normal): in each of these two cases, there is only one
isomorphism type of groups of this kind. Both groups satisfy condition (i) of
Theorem 1.2.

 It will be convenient to have a short temporary name for the $2$-groups $G$
such that if $g,h\in G$ then $\gp{g^2}$ is normal in $G$ and $\gp{g,h}/\gp{g^2}$
is either abelian or dihedral: let us call these groups $G$ {\it good.}
Obviously, a group is good if and only if each of its two-generator subgroups is
good, and so all subgroups of good groups are good. A little more thought shows
that all factor groups of good groups are also good, and that the direct product
of an elementary abelian $2$-group with a good group is always good.

 Of course, all abelian or dihedral $2$-groups are good. The next exploratory
step is to look (for example, by using \cite{\mhs}) at each of the groups of
order dividing $16$, and check that all but three of them are good. The three
that  fail do so because they are of the form $\gp{g,h}$ with $g^2=1$ but are
neither abelian nor dihedral; they are $(C_2\times C_2)\rtimes C_4$ and the two
semidirect products $C_8\rtimes C_2$ in which the action of $C_2$ on $C_8$ is
neither trivial nor inverting. (The nonabelian $(C_2\times C_2)\rtimes C_4$ form
a single isomorphism class.) We note for future use that the reasons which  make
the generalized quaternion group of order $16$ good but the semidihedral group
 of
order $16$ bad, yield the same conclusions for generalized quaternion groups and
semidihedral groups of larger orders as well.

\smallskip

\noindent{\it Proof of the last sentence of Lemma \rom{1.4}.\quad}In case (i),
 all
two-generator subgroups of $H$ are abelian or dihedral so $G$ is good. In cases
(ii)-(v) we have $|\Phi(H)|\leq4$, so the two-generator subgroups $K$ of $H$ are
of order dividing $16$. No $K$ can be a $C_8\rtimes C_2$, because that would
 mean
$\Phi(H)\geq\Phi(K)\cong C_4$ but $\Phi(H)$ contains no $C_4$.  No $K$ can be a
$(C_2\times C_2)\rtimes C_4$, because then
$\Phi(H)\geq\Phi(K)\cong C_2\times C_2$ would exclude case (ii) and so (see
 Lemma
4.1) ensure $|\Omega(H)|=4$, contrary to $\Omega(K)\leq\Omega(H)$ and
$|\Omega(K)|=8$. Thus in every case we may conclude that $G$ is good. This
 proves
the converse part of Lemma 1.4. \qed

\smallskip

 The proof of the direct part needs more preparation. To avoid cumbersome
circumlocution, {\it we count $Q_8$ among the generalized quaternion groups.} As
usual, $x^y=y\m1xy$, and an {\it involution} is a group element of order $2$.

 The first step is rather trivial, and the second is not much harder.

\lemma{4.2}
 {If $H$ is a nonabelian good group with $|\Phi(H)|=2$ and no direct factor of
order $2$, then it satisfies \rom{(ii)} of Theorem \rom{1.2}.}

\proof
 Since $H$ has no direct factor of order $2$, it has no central involution
outside $\Phi(H)$; thus $|\Phi(H)|=2$ implies that $\sf Z(H)$ is cyclic of order
at most $4$. If $|\sf Z(H)|=2$ then $H$ is extraspecial; otherwise there is a
maximal subgroup $M$ which does not contain $\sf Z(H)$ and is easily seen to be
extraspecial.\qed

\lemma{4.3}
 {An involution in a good group normalizes every cyclic subgroup of order
 greater
than $2$ and centralizes every subgroup isomorphic to $C_4\rtimes C_4$.}

\proof
 Let $g$ be an involution in a good group $G$. If $h$ is any element of $G$,
then by the definition of `good' $\gp{g,h}$ is abelian or dihedral: in either
case, if $\gp{h}$ is of order greater than $2$ then it is normalized by $g$.
Suppose now that $x$, $y$ are elements of $G$ such that
$\gp{x,y}\cong C_4\rtimes C_4$ with $x^y=x\m1$. If $x^g=x\m1$ and
$y^g=y^{\pm1}$, then $g$ fails to normalize $\gp{xy}$, while if $x^g=x$ and
$y^g=y\m1$, then $yg$ is an involution which does not normalize $\gp{xy}$: so
the only option is that $g$ centralizes $\gp{x,y}$. \qed

\smallskip

 We shall make repeated use of the fact that if $\gp{x,y}$ is a nonabelian
dihedral $2$-group then $\gp{x}$, $\gp{y}$, $\gp{xy}$ are pairwise distinct and
precisely two of them are non-normal subgroups of order $2$, while the third is
normal and has order divisible by $4$.

 If $\gp{g,h}$ is good and $\gp{g}$ is not normal in it, then
$\gp{g,h}/\gp{g^2}$ is a nonabelian dihedral $2$-group, so it follows that
$\gp{g^2}$ contains precisely one of $\gp{h^2}$ and $\gp{(gh)^2}$. Suppose
further that neither $\gp{h}$ nor $\gp{gh}$ is normal in $\gp{g,h}$: then, by
this argument, each of $\gp{g^2}$, $\gp{h^2}$ and $\gp{(gh)^2}$ must contain one
and only one of the other two. As it is impossible to order a three-element set
in this manner, we have a contradiction, which proves that at least one of
$\gp{g^2}$, $\gp{h^2}$ and $\gp{(gh)^2}$ must be normal in $\gp{g,h}$. We have
proved that {\it in a good group, every two-generator subgroup is metacyclic}
(in the sense of having a cyclic normal subgroup with cyclic quotient).

 If a nonabelian $2$-group has a cyclic normal subgroup of order $4$ with cyclic
quotient, then it is isomorphic to one of
$$
\eqalign{
 P_k&=\bigl\langle\,w,y\bigm|w^4=1,\ w^y=w\m1,\ y^{2^k}=w^2\,\bigr\rangle,\cr
 R_k&=\bigl\langle\,w,y\bigm|w^4=1,\ w^y=w\m1,\ y^{2^{k+1}}\!=1\,\bigr\rangle.
\cr}
$$
 (This notation is not intended for use beyond the proof of the next lemma.) We
claim that {\it such a group is good if and only if} $k\leq1$. If $k\leq1$ then
both groups are of order dividing $16$ and we have said nothing new. If $k>1$
then $P_k=\bigl\langle wy^{2^{k-1}},y\bigr\rangle$ and
$\bigl(wy^{2^{k-1}}\bigr)^2=1$ but $P_k$ is neither abelian nor dihedral and  so
cannot be good. As $P_k$ is a homomorphic image of $R_k$, in this case $R_k$
cannot be good either.

\lemma{4.4}
 {If $\gp{x,y}$ is good and $\gp{x}$ is normal in it but $\gp{y}$ is not, then
$|x|\geq4\geq|y|$ and $x^y=x\m1$.}

\proof
 Since in $\gp{x,y}/\gp{y^2}$ the image of $\gp{x}$ is normal but the image of
$\gp{y}$ is not, $\gp{x,y}/\gp{y^2}$ is a nonabelian dihedral group and the
 image
of $\gp{x}$ has order divisible by $4$: thus
$$
 x^2\notin\gp{y^2}\eqno(1)
$$
 and
$$
 x^y\equiv x\m1\ \bmod\gp{y^2}.\eqno(2)
$$
 Of course (2) and the normality of $\gp{x}$ give that
$$
 x^y\equiv x\m1\ \bmod\ \gp{x}\cap\gp{y^2}\,.\eqno(3)
$$
 From (1) we know that $\gp{x}\cap\gp{y^2}\leq\gp{x^4}$, so by (3)
$$
 x^y\equiv x\m1\ \bmod\gp{x^4}\,.\eqno(4)
$$
 It also follows from (1) that there is in $\gp{x}$ and element, $w$ say, of
order $4$, and (4) implies that $w^y=w\m1$. Thus $\gp{w,y}$ is a $P_k$ or an
$R_k$, and the argument leading up to the lemma may be invoked for the
 conclusion
that $y^4=1$. If $y^2\notin\gp{x}$, then
$\gp{x}\cap\gp{y^2}=1$ and so (3) gives that $x^y=x\m1$. If
$y^2\in\gp{x}$, then $\gp{x,y}$ has a cyclic maximal subgroup and a nonabelian
dihedral quotient, so it can only be dihedral or semidihedral or generalized
quaternion. We have already seen that semidihedral groups are not good, so
$x^y=x\m1$ holds in this case as well. Finally, $|x|\geq4$ because $\gp{x,y}$ is
nonabelian. This completes the proof  of the lemma.\qed

\smallskip

It follows that under the hypotheses of Lemma 4.4 the group $\gp{x,y}$ is either
dihedral or generalized quaternion or a semidirect product
$\gp{x}\rtimes\gp{y}$ with $|x|\geq|y|=4$ and $x^y=x\m1$.

\lemma{4.5}
 {f $\gp{x,u}$ is good and both $\gp{x}$ and $\gp{u}$ are normal in it, then
$\gp{x,u}$ is either abelian or isomorphic to $Q_8$.}

\proof
 We argue by contradiction. If the conclusion is false then $\gp{x,u}$ cannot be
hamiltonian: so there is a cyclic subgroup, $\gp{y}$ say, which at least one of
$x$ and $u$ fails to normalize. We may assume without loss of generality that
$x$ does not normalize $\gp{y}$, and then Lemma 4.4 is conveniently applicable.
At first we only use the conclusion that $x^y=x\m1$. By conjugation, $\gp{x,u}$
induces a cyclic group of automorphisms on $\gp{x}$, and now we know that this
includes the inverting automorphism. In the automorphism group of a cyclic
$2$-group, the subgroup generated by the inverting automorphism is a maximal
cyclic subgroup ($-1$ is not a square \,mod~$2^n$ when $n>1$): so the group
induced by $\gp{x,u}$ is of order $2$. This proves that the centralizer
$\sf C_{\gp{x,u}}(x)$, which contains $x$ but not $y$, is of index $2$. Thus
$y\notin\gp{x,\Phi(\gp{x,u})}$, and therefore $\gp{x,y}=\gp{x,u}$. However, of
the groups of Lemma 4.4, only $Q_8$ can be generated by two normal cyclic
subgroups, and we have assumed that $\gp{x,u}$ is not that group. This
contradiction completes the proof of Lemma 4.5.\qed

\smallskip

 Since a good two-generator group is metacyclic, it has a generating set which
satisfies the hypotheses of one of these two lemmas.

\corollary{4.6}
 {A two-generator $2$-group is good if and only if it is either abelian or
dihedral or generalized quaternion or a semidirect product $\gp{x}\rtimes\gp{y}$
with $|x|\geq|y|=4$ and $x^y=x\m1$.\qed}

\lemma{4.7}
 {If $H$ is a nonabelian good group of exponent greater than $4$, then $H$
satisfies \rom{(i)} of Theorem \rom{1.2}.}

\proof
 In the semidirect products of Corollary 4.6, all the elements outside the
abelian group $\gp{x,y^2}$ have order $4$, and every cyclic subgroup of
$\gp{x,y^2}$ is normal in $\gp{x,y}$. In a dihedral or generalized quaternion
group which does not have exponent $4$, all cyclic subgroups of order greater
than $4$ are normal and lie in the unique cyclic maximal subgroup. It follows
that in a good group every cyclic subgroup of order greater than $4$ is normal
and any two elements of order greater than $4$ commute.

 Let $A=\gp{\,a\in H\mid a^4\neq1\,}$: this is now an abelian subgroup of $H$.
Let $a\in H$ with $a^4\neq1$, and $h\in H$ but $h\notin A$, so $h^4=1$. Then $a$
and $h$ cannot commute (else $(ah)^4\neq1$ and hence $ah\in A$, $h\in A$ would
follow). Lemma~4.5 cannot apply with $x=a$, $u=h$, because $\gp{a,h}$ is neither
abelian nor of exponent $4$. Hence Lemma 4.4 must apply with $x=a$, $y=h$. It
follows that every element of $H$ outside $A$ must invert every element of $A$.
If $H$ had more than one nontrivial coset modulo $A$, the quotient of two
elements chosen from distinct nontrivial cosets would still lie outside $A$: it
would have to fix as well as invert every element of $A$. This being impossible,
the index of $A$ in $H$ must be $2$. \qed

\lemma{4.8}
 {If $H$ is a good group of exponent $4$, then its Frattini subgroup is
 elementary
abelian and central.}

\proof
 If $g,h\in H$ then $[g^2,h]=1$ because, by Corollary 4.6, $\gp{g,h}$ is either
abelian or $D_8$ or $Q_8$ or $C_4\rtimes C_4$, and $[g^2,h]=1$ holds for every
pair of elements $g,h$ in each of these groups. This shows that the Frattini
subgroup is generated by central involutions. \qed

\smallskip

 In the proof of our next lemma, we shall make use of two properties of
$C_4\rtimes C_4$. First, it has only two nontrivial elements that are squares.
Second, as it can be generated by two non-commuting elements of order
$4$ whose product is also of order $4$, no automorphism of it can invert all
elements of order $4$.

\lemma{4.9}
 {If $a,x,y$ are elements of a good group $H$ of exponent $4$ such that
$x^y=x\m1$ and $a^2\notin\gp{x,y}\cong C_4\rtimes C_4$, then $x$ must centralize
and $y$ must invert $a$. If also $b\in H$ and $a^2\neq b^2\notin\gp{x,y}$, then
$a$ and $b$ commute.}

\proof
 The subgroup $\gp{a^2}$ is normal and the image of $\gp{x,y}$ in the quotient
$H/\gp{a^2}$ is still a $C_4\rtimes C_4$. By Lemma 4.3, this image must
centralize the image of $a$. It follows that $\gp{x,y}$ normalizes $\gp{a}$. It
cannot centralize $a$, for then we would have
$\gp{a,x,y}=\gp{a}\times\gp{x,y}$, and we know that
$C_4\times(C_4\rtimes C_4)$ has a quotient $C_4\times D_8$ which is not good.
Since $\gp{xy,y}=\gp{x,y}$, it follows that at least one of $xy$ and $y$ must
invert $a$.

 Suppose only one of them does: say, $a^{xy}=a\m1$ but $a^y=a$. Then
$(ay)^2\notin\gp{x,y}$, and the above argument may be repeated with $ay$ in
place of $a$, giving the conclusion that at least one of $xy$ and $y$ must
invert $ay$. However, now $(ay)^{xy}\neq(ay)\m1$ because
$a\m1x^2y\neq a\m1y\m1$, and $(ay)^{y}=ay\neq(ay)\m1$: we have reached a
contradiction. A similar argument gives a contradiction if we assume that
$a^{xy}=a$ and $a^y=a\m1$. This proves that both $xy$ and $y$ must invert $a$,
that is, $x$ must centralize and $y$ must invert $a$.

 By Corollary 4.6, the only nonabelian good groups of exponent $4$ generated by
two elements of order $4$ are $Q_8$ and $C_4\rtimes C_4$. Since $a^2\neq b^2$,
we cannot have $\gp{a,b}\cong Q_8$. If $\gp{a,b}\cong C_4\rtimes C_4$, then
$a^2$ and $b^2$ are the only nontrivial squares in $\gp{a,b}$, and by assumption
neither of these lies in $\gp{x,y}$: thus all cyclic subgroups of order $4$ in
$\gp{a,b}$ avoid $\gp{x,y}$ and are therefore inverted by $y$. Since no
automorphism of $C_4\rtimes C_4$ can act like that, $a$ and $b$ must
 commute.\qed

\smallskip

 One of the two nontrivial squares in $C_4\rtimes C_4$ generates the commutator
subgroup; hence if two cyclic subgroups of order $4$ in $C_4\rtimes C_4$
intersect trivially, one of them must be normal. This will also be used in the
proof of the next lemma.

 In view of Lemma 4.8, if $H$ is a good group of exponent $4$ then $\Phi(H)$ is
an elementary abelian group spanned by squares, so it has a basis consisting of
squares: that is, $H$ has a subset $X$ such that
$\Phi(H)=\prod_{x\in X}\gp{x^2}$ and each $x^2$ is nontrivial.

\lemma{4.10}
 {Let $H$ be a good group of exponent $4$ and $X$ any subset such that
$\Phi(H)=\prod_{x\in X}\gp{x^2}$ and each $x^2$ is nontrivial. Then either
$\gp{X}$ is abelian or all but one of the elements of $X$ commute with each
other and are inverted by the remaining one.}

\proof
 Suppose that $\gp{X}$ is nonabelian, and that $x,y$ is a noncommuting pair of
elements of $X$. Since $\gp{x}\cap\gp{y}=1$, one of $\gp{x}$ and $\gp{y}$
normalizes the other: say, $x^y=x\m1$. By Lemma 4.9, then each element of
$X\setminus\{x,y\}$ is centralized by $x$ and inverted by $y$, and any two
elements of $X\setminus\{x,y\}$ commute with each other. \qed

\lemma{4.11}
 {If $H$ is a nonabelian good group of exponent $4$, then the Frattini subgroup
of its centre has order at most $2$.}

\proof
 Suppose not: then $\sf Z(H)=\gp{a}\times\gp{b}\times\cdots$ with $|a|=|b|=4$.
Since $H$ is nonabelian, it has a nonabelian two-generator subgroup $K$. By
Corollary~4.6, $K$ is either $D_8$ or $Q_8$ or $C_4\rtimes C_4$. We propose to
show that $\gp{a,b,K}$ must contain a subgroup with a quotient isomorphic to
$C_4\times D_8$. Since $C_4\times D_8$ has a subgroup
$(C_2\times C_2)\rtimes C_4$ which we know is not good, this will contradict the
assumption that $H$ is good, and so prove the lemma. If $\gp{a^2,b^2}\leq K$, we
must have $K\cong C_4\rtimes C_4$, and then no generality is lost by assuming
that $K=\gp{\,x,y\mid x^4=y^4=1,\ x^y=x\m1\,}$ and $a^2=x^2$, $b^2=y^2$, so
$\gp{a,ax,by}=\gp{a}\times\gp{ax,by}=C_4\times D_8$. This argument tacitly
involved changing $a$ and $b$ if necessary (without changing $\gp{a,b}$). The
same flexibility allows us to assume that if $\gp{a^2,b^2}\nleq K$ then
$\gp{a,b,K}=\gp{a}\times\gp{b,K}$, and then what we need is that $D_8$ is a
quotient of a subgroup of $\gp{b,K}$. Since $D_8$ is a quotient of
$C_4\rtimes C_4$, this is only an issue if $K\cong Q_8$, but then the central
product $C_4\sf Y\,Q_8$ is a quotient of $\gp{b,K}$ and of course
$C_4\sf Y\,Q_8\cong C_4\sf Y\,D_8>D_8$. \qed

\lemma{4.12}
 {If $H$ is a nonabelian good group of exponent $4$ with $|\Phi(H)|>4$, then $H$
satisfies \rom{(i)} of Theorem \rom{1.2}.}

\proof
 Let $X$ be a subset of $H$ of the kind discussed in Lemma 4.10.

 First, suppose that $\gp{X}$ is abelian; then $\gp{X}=\prod_{x\in X}\gp{x}$.
Set $A=\sf C(\gp{X})$: by Lemma 4.11, this centralizer is abelian. As $A$ is of
exponent precisely $4$, it is generated by its elements of order $4$. Thus if an
element $h$ normalizes every cyclic subgroup of order $4$ in $A$ then it either
centralizes or inverts $A$. If $h$ is an involution, then by Lemma 4.3 this
comment is applicable. If we can show that each $h$ of order $4$ acts on
$A$ in this way, the claim of the lemma will follow.

 Suppose then that $h\in H$, $h\notin A$, and $|h|=4$. If necessary, one can
change $X$ without changing $\gp{X}$ (and therefore without changing $A$) so  as
to achieve that $X$ has an element, $x_1$ say, with $x_1^2=h^2$. Let $x_2$ be
another element of $X$, and set $X'=X\setminus\{x_1,x_2\}$. Then
$\{h,x_2\}\cup X'$ and $\{h,x_1x_2\}\cup X'$ can also play the role of $X$ in
Lemma 4.10. Since $X'$ is nonempty (this is where we use the assumption that
$|\Phi(H)|>4$) and commutes with $x_2$ and with $x_1x_2$, in each case $h$ is
the only element which could invert all the others. Since $X'$ is a nonempty
common part of `all the others', $h$ behaves the same way in both cases. If it
centralizes in both cases, then it centralizes all of $X$, contrary to
$h\notin A$. Thus $h$ inverts all elements of $\gp{X}$.

 This proves that the centralizer $A$ of $X$ has at most one nontrivial coset in
$H$, so $|H:A|$ is at most $2$. It cannot be $1$, because $A$ is abelian but
$H$ is not. If $a$ is any element of order $4$ in $A$, then there is an $x$ in
$X$ such that for $Y=\{a\}\cup(X\setminus\{x\})$ we have
$\Phi(H)=\prod_{y\in Y}\gp{y^2}$, so $Y$ can play the role of $X$ in all this.
The centralizer of $Y$ contains and therefore equals $A$, so an element $h$ of
order $4$ outside $A$ inverts every element of $\gp{Y}$ as well. This proves
that such an $h$ inverts every element of order $4$ in $A$, and so it inverts
every element of $A$, as required.

 Second, suppose that $\gp{X}$ is not abelian: say, $y$ is the element of
$X$ which inverts all the others, and all the others commute with each other.
Set $X'=X\setminus\{y\}$ and $A=\sf C(X')$. By Lemma 4.11, $A$ is abelian. If
$\Phi(A)=\Phi(H)$, we can replace $X$ by a subset of $A$ and appeal to the half
of the lemma which we have already proved. It remains to deal with the case
$\Phi(A)<\Phi(H)$. Of course then $\Phi(A)=\Phi(\gp{X'})$, so $A$ is the direct
product of $\gp{X'}$ with an elementary abelian group. We shall show that every
element $h$ outside $A$ inverts $A$.

 Now $\gp{x,y}\cong C_4\rtimes C_4$ whenever $x\in X'$, so Lemma 4.3 ensures
that every involution lies in $A$: we need only consider the $h$ of order $4$.
If $h^2\notin\bigl\langle\,x^2\bigm|x\in X'\,\bigr\rangle$ then $\{h\}\cup X'$
can play the role of $X$ in Lemma 4.11; given that $X'$ has at least two
elements and they commute, this means that $h$ must either centralize or invert
every element of $X'$, that is, $h\in\gp{y}A$. If
$h^2\in\bigl\langle\,x^2\bigm|x\in X'\,\bigr\rangle$, one can change $X'$
without changing $\gp{X'}$ (and therefore without changing $A$) so  as to
achieve that $X'$ has an element, $x_1$ say, with $x_1^2=h^2$. Let $x_2$ be
another element of $X'$ (here we use again the assumption that
$|\Phi(H)|>4$), and set $X''=X'\setminus\{x_1,x_2\}$. Both
$\{y,h,x_2\}\cup X''$ and $\{y,h,x_1x_2\}\cup X''$ can play the role of $X$ in
Lemma 4.11. Because $y$ inverts both $x_2$ and $x_1x_2$, we can conclude that
$h$ commutes with $x_2$, with $x_1x_2$, and with every element of $X''$ (and is
inverted by $y$). What matters is that in this case $h\in A$.

 We have proved that $H=\gp{y}A$, and noted that $y$ centralizes, that is,
inverts, every involution. We have also seen that $y$ inverts $\gp{X'}$. Since
$A$ is the direct product of $\gp{X'}$ with an elementary abelian group, it
follows that $y$ inverts $A$, and then so does every element of $H$ outside $A$.
This completes the proof of the lemma.\qed

\lemma{4.13}
 {If $h$ is a noncentral involution in a good group $H$ of exponent $4$ with
$|\Phi(H)|>2$, then $H$ has an abelian subgroup $A$ of index $2$ such that every
element of $A$ is inverted by $h$.}

\proof
 Since $h$ is noncentral it is noncentral already in some nonabelian
two-generator subgroup which by Corollary 4.6 can only be dihedral: thus there
is in $H$ an element $a$ of order $4$ such that $a^h=a\m1$. Since
$|\Phi(H)|>2$, there is also in $H$ an element $b$ such that
$a^2\neq b^2\neq1$. For any such $b$, by Corollary~4.6, $\gp{a,b}$ is either
abelian of a $C_4\rtimes C_4$. It cannot be the latter, because $h$ does not
centralize it and we have Lemma 4.3. It follows that every such $b$ commutes
with $a$. Further by Lemma 4.3, $h$ must normalize both $\gp{b}$ and $\gp{ab}$,
which is now only possible if $h$ inverts $b$. If $c\in H$ and $c^2\neq1$, then
either $c^2\neq a^2$ or $c^2\neq b^2$, and the above argument with $a,c$ or
$b,c$ in place of $a,b$ yields that $h$ inverts $c$. We have proved that $h$
inverts every element of order $4$ in $H$. Further, any two elements of order
$4$ commute: else by Corollary 4.6 the subgroup they generate would be a $Q_8$
or a $C_4\rtimes C_4$, and we have observed just before stating this lemma that
neither of these groups has an automorphism that inverts all elements of order
$4$. Set $A=\gp{\,a\in H\mid a^2\neq1\,}$; this is an abelian subgroup, and
every element of it is inverted by $h$. If $g$ is an element of $H$ outside
$A$, then $g$ is an involution (by the definition of $A$). If $\gp{a,g}$ were
abelian, it would be generated by elements of order $4$ and so would lie in
$A$, contrary to $g\notin A$: thus $g$ is a noncentral involution and, like
$h$, inverts every element of $A$. Hence $gh$ centralizes $A$ and therefore
cannot lie outside it. This means that every element $g$ of $H$ outside $A$ lies
in the coset $Ah$, that is, that $|H:A|$=2. \qed

\smallskip

\noindent{\it Proof of the direct part of Lemma \rom{4.1}.\quad}Let $G$ be a
nonabelian good group. By Lemma 2.3, we may write $G$ as $E\times H$ with $E$
elementary abelian and $H$ having no direct factor of order $2$, and of course
$H$ is also nonabelian and good. If the exponent of $H$ is greater than $4$,
Lemma 4.7 shows that $H$ satisfies (i) of Theorem 1.2. Suppose that the exponent
of $H$ is $4$. By Lemma 4.8, then $\Phi(H)\leq\sf Z(H)$. If $|\Phi(H)|=2$ then
Lemma 4.2 shows that $H$ satisfies (ii), while if $|\Phi(H)|>4$ then (i) holds
 by
Lemma 4.12. In the remaining case, $|\Phi(H)|=4$. If $H$ has a noncentral
involution, (i) holds by Lemma 4.13. If all involutions are central, then
$\Omega(H)\leq\Phi(H)$ because $H$ has no direct factor of order $2$. We cannot
have $|\Omega(H)|=2$, for then $H$ would be cyclic or generalized quaternion and
(as $H$ has exponent $4$) this is excluded by $|\Phi(H)|=4$. Thus
$\Omega(H)=\Phi(H)$. Let $\Phi(H)=\gp{g^2,h^2}$ with $g,h\in H$, and $K$ any
subgroup of $H$ which contains $\gp{g,h}$: then also
$\gp{g^2,h^2}=\Omega(K)=\Phi(K)\leq\sf Z(K)$ and so Lemma 4.1 (applied to $K$ in
place of $H$) shows that $|K|\neq128$. It follows that $H/\gp{g,h}$ is an
elementary abelian $2$-group without a subgroup of order $128/16$; hence $H$ has
order dividing $64$. Lemma 4.1 therefore shows that $H$ satisfies (i) or (iii)
 or
(iv) or (v).

 This completes the proof of Lemma 1.4, and so also the proof of Theorem~1.2.
\qed

\medskip

 The application of Lemma 4.1 to $K$ in place of $H$ raises a question: had the
inspection of that CAYLEY library produced a different answer, would we still
have a theorem? If the list of Lemma 4.1 did contain at least one group of order
$128$, we could not save the present proof simply by adjusting the list in
Theorem 1.2. However, instead of applying Lemma 4.1 to $K$ we could in any case
appeal to Lemma 4.14 below, and so deduce that in the remaining case $H$ has
order dividing $128$. After that, even a Lemma 4.1 with a modified list would be
good enough to conclude the proof (of a suitably modified theorem).

\lemma{4.14}
  {If $H$ is a finite $2$-group such that $\Phi(H)\leq\Omega(H)\leq\sf Z(H)$ and
$|\Omega(H)|\leq2^n$, then $|H|\leq2^{n(n+5)/2}$.}

\proof
 If $n=0$ then $H=1$ while if $n=1$ then $H\leq Q_8$, so we have the initial
step for a proof by induction on $n$. For the inductive step, suppose that
$n>1$. If $H$ has exponent $2$, the claim is obvious. Let $h$ be an element
of order $4$ in $H$. Since the conjugates of $h$ differ from $h$ by commutators
and there are at most $2^n$ commutators in $H$, the centralizer $\sf C(h)$ has
index at most $2^n$. As $\Phi(H)$ is of exponent $2$, it cannot contain $h$. Let
$M$ be a maximal subgroup which does not contain $h$. Set $K=M\cap\sf C(h)$, and
note that $|H:K|\leq2^{n+1}$. If $h^2=k^2$ for some element $k$ of $K$, then
$(hk\m1)^2=1$ and so $hk\m1\in\Omega(H)=\Phi(H)\leq M$ contradicts $h\notin M$.
Thus all involutions of $K/\gp{h^2}$ lie in the group $\Omega(K)/\gp{h^2}$, and
then by the inductive hypothesis $|K/\gp{h^2}|\leq2^{(n-1)(n-1+5)/2}$. Together
with $|H:K|\leq2^{n+1}$, this inequality gives the bound on $|H|$ that we wanted
and so completes the inductive step.\qed

\smallskip

 At $n=2$, from Lemma 4.1 we can get $2^6$ where Lemma 4.14 gives $2^7$. Running
the proof of Lemma 4.14 with this starting point, we can get an improved bound
for higher values of $n$ as well. On the other hand, the direct product of $n$
copies of $Q_8$ shows that the $n(n+5)/2$ in  Lemma 4.14 could  not be lowered
below $3n$. One may well wonder just what the optimal bound is.

\begref{References}
\refno {\aab.}A.\,A. Bovdi,
Unitarity of the multiplicative group of an integral group ring.
Mat. Sb. {\bf 119}, 384--400 (1982) (Russian);
English transl. in Math. USSR Sbornik {\bf 47}, 377--389 (1984)
\refno {\bgs.}A.\,A. Bovdi, P.\,M. Gudivok, M.\,S. Semirot,
Normal group rings.
Ukrain. Mat. Zh. {\bf 37}, 3--8 (1985) (Russian)
\refno {\bsb.}A.\,A. Bovdi and S.\,K. Sehgal,
Unitary subgroups of integral group rings.
Publ. Mat. {\bf 36}, 197--204 (1992)
\refno {\bsn.}A.\,A. Bovdi and S.\,K. Sehgal,
Unitary subgroups of integral group rings.
Manuscripta Math. {\bf 76}, 213--222 (1992)
\refno {\vab.}V.\,A. Bovdi,
Normal twisted group rings.
Dokl. Akad. Nauk Ukrain. SSR (1990) No 7, \ 6--7 (Russian)
\refno {\mhs.}Marshall Hall, Jr. and James K. Senior,
The groups of order $2^n$ $(n\leq6)$.
Macmillan, New York, 1964
\refno {\hig.}Graham Higman,
Suzuki $2$-groups.
Illinois J. Math. {\bf 7}, 79--96 (1963)
\refno {\kap.}Irving Kaplansky,
Infinite abelian groups.
Revised edition, second printing;
University of Michigan Press, Ann Arbor, 1971
\refno {\mfn.}M.\,F. Newman,
On a class of nilpotent groups.
Proc. London Math. Soc. $(3)$ {\bf 10}, 365--375 (1960)
\refno {\nob.}M.\,F. Newman and E.\,A. O'Brien,
A CAYLEY library for the groups of order dividing $128$.
Proceedings of the Singapore Group Theory Conference held
at the National University of Singapore, 1987,
ed. by Kai Nah {\smc Cheng} and Yu Kiang {\smc Leong},
de Gruyter, Berlin, New York, 1989; 437--442
\refno {\rse.}J\"urgen Ritter and Sudarshan K. Sehgal,
Generators of subgroups of $U(ZG)$.
Contemp. Math. {\bf 93}, 331--347 (1989)
\endref

\address{Victor Bovdi
Department of Mathematics
Bessenyei Teachers' College
4401 Ny\'\i regyh\'aza
Hungary}

\address{L.\,G. Kov\'acs
School of Mathematical Sciences
Australian National University
Canberra ACT 0200
Australia}

\bye